\documentclass[a4paper,english,oneside]{article}
\usepackage[T1]{fontenc}
\usepackage[utf8]{inputenc} 

\usepackage{geometry}
\usepackage{amsmath}
\usepackage{amsfonts}
\usepackage{amssymb}
\usepackage{amsthm}
\usepackage[english]{babel}
\usepackage{color}
\usepackage{lmodern}
\usepackage{pstricks-add}
\usepackage{hyperref}
\hypersetup{hidelinks}

\usepackage{mathscinet}
\usepackage{enumitem}
\usepackage[cal=boondoxo,bb=ams]{mathalfa}
\usepackage{microtype}

\usepackage{mathtools}
\usepackage{esint}

\theoremstyle{plain}
\newtheorem{theorem}{Theorem}[section]
\newtheorem{lemma}[theorem]{Lemma}

\newtheorem{proposition}[theorem]{Proposition}

\theoremstyle{definition}
\newtheorem{definition}[theorem]{Definition}

\theoremstyle{remark}
\newtheorem{remark}{Remark}

\begin{document}

\title{Semilinear damped wave equation on a compact Lie group \\ with a non-autonomous forcing term}

\author{Wenhui Chen$^\mathrm{a}$, Sandra Lucente$^\mathrm{b}$, Alessandro Palmieri$^{\mathrm{c}}$\footnote{email addresses: wenhui.chen.math@gmail.com, sandra.lucente@uniba.it, alessandro.palmieri@uniba.it}} 

\date{}
 \date{$^\mathrm{a}$\small{School of Mathematics and Information Science, Guangzhou University, 510006 Guangzhou, P.R. China}\\
 $^\mathrm{b}$\small{Department of Physics, University of Bari, 70125 Bari, Italy}\\
 $^\mathrm{c}$\small{Department of Mathematics, University of Bari, 70125 Bari, Italy}}

\maketitle

\begin{abstract}

In the present note, we consider a semilinear damped wave equation on a compact Lie group with a non-autonomous nonlinearity $\varphi(t)|u|^p$. We are interested in describing how the nonnegative time-dependent factor $\varphi$ effects the global in time prolongability of a local solution. In particular, the summability of the function $\varphi$ provides a criterion to distinguish between the blow-up in finite time and global existence of small data. Finally, we derive sharp lifespan estimates for local in time solutions when $\varphi \not\in L^1([0,+\infty))$ and satisfies a certain scaling condition, that we named uniform upper scaling condition.
\end{abstract}

\begin{flushleft}
\end{flushleft}

\begin{flushleft}
\textbf{AMS Classification (2020)}:  35B44, 35B60, 35L71 
\end{flushleft}

\section{Introduction}

Let $\mathbb{G}$ be a compact Lie group and let $\mathcal{L}$ be the Laplace-Beltrami operator on $\mathbb{G}$.
In this paper, we investigate blow-up results and global existence results for small data solutions for the following Cauchy problem 
\begin{align}\label{semilinear CP damped phi}
\begin{cases} \partial_t^2 u-\mathcal{L} u+\partial_t u =\varphi(t)| u|^p, &  x\in \mathbb{G}, \ t>0,\\
u(0,x)=\varepsilon u_0(x), & x\in \mathbb{G}, \\ \partial_t u(0,x)=\varepsilon u_1(x), & x\in \mathbb{G},
\end{cases}
\end{align} where $p>1$, $\varepsilon$ is a positive constant describing the size of Cauchy data and $\varphi \in L^1_{\mathrm{loc}}([0,+\infty))$ is a nonnegative function. Our main goal is to study how the time-dependent factor $\varphi$ in the non-autonomous right-hand side does influence the global in time prolongability of local in time solutions to the semilinear damped wave equation in \eqref{semilinear CP damped phi}. We emphasize that the summability property of $\varphi$ (respectively its lack of) is crucial to guarantee the existence of global in time small data solutions to \eqref{semilinear CP damped phi} (respectively the blow-up in finite time). We are going to prove two types of results: on the one hand, local/global existence results for energy solutions obtained by employing $L^2(\mathbb{G})-L^2(\mathbb{G})$ estimates for the corresponding linear homogeneous problem in a standard contraction argument; on the other hand,  a blow-up result is proved by using the classical approach for semilinear wave models (with nonnegative power nonlinearities) in the Euclidean framework for a not summable $\varphi$ satisfying a \emph{uniform upper scaling condition} (cf. Definition \ref{Def UUSC}). In particular, in the case $\varphi\not \in L^1([0,+\infty))$ we shall obtain the same lower bound estimates (in the local existence result) and upper bound estimates (in the blow-up result) for the lifespan of a solution blowing up in finite time up to different multiplicative constants independent of $\varepsilon$.

Let us begin by reviewing the literature for the case with autonomous nonlinear term ($\varphi\equiv 1$). For the classical semilinear damped wave equation in $\mathbb{R}^n$ 
\begin{align*}
\begin{cases} \partial_t^2 u-\Delta u+\partial_t u =| u|^p, &  x\in \mathbb{R}^n, \ t>0,\\
u(0,x)=\varepsilon u_0(x), & x\in \mathbb{R}^n, \\ \partial_t u(0,x)=\varepsilon u_1(x), & x\in \mathbb{R}^n,
\end{cases}
\end{align*}
the critical exponent is the so-called Fujita exponent $p_{\mathrm{Fuj}}(n)\doteq 1+\frac{2}{n}$ , cf. \cite{Mat76,TY01,Zhang01,IT05}. This exponent coincides with the critical exponent for the Cauchy problem associated to the semilinear heat equation $\partial_t u-\Delta u=|u|^p$ and it is named after the author of \cite{Fuj66}. The global existence of small data solutions in the supercritical case is proved for $n=1,2$ in \cite{Mat76}. Afterwards, the global existence  was proved for any $n\geqslant 1$ in the supercritical case $p>p_{\mathrm{Fuj}}(n)$ by working with compactly supported initial data in \cite{TY01} and later in \cite{IT05} without compact supports for the data. Moreover, in \cite{TY01}  the blow-up of local in time solutions is proved (under suitable sign assumptions for the Cauchy data) in the subcritical case $1<p<p_{\mathrm{Fuj}}(n)$. Then, in \cite{Zhang01} the blow-up is proved in the critical case $p=p_{\mathrm{Fuj}}(n)$ as well. 

The Fujita-type exponent has been also found as critical exponent for the semilinear heat equation on stratified Lie groups (see \cite{Pas98} and \cite{GP19} for the sharp lifespan estimates in the Heisenberg group) and on more general unimodular Lie groups in \cite{RY18}. In particular, on compact Lie groups the critical exponent is $\displaystyle{\lim_{D\to 0^+} p_{\mathrm{Fuj}}(D)=+\infty}$.

As in the Euclidean case, these critical exponents for the semilinear heat equations are the same for the corresponding semilinear damped wave equations on the Heisenberg group \cite{Pal19,GP19DW} and on a compact Lie group \cite{Pal20DWE}.
Recently, other papers have been devoted to semilinear evolution models on compact Lie groups (see, for instance, \cite{Pal20WE,Pal20WEdm,BKM24,ChenPal26}).

Let us consider now the case with a non-autonomous forcing term.
In this case, fewer results are known in the literature.

Concerning the semilinear heat equation $\partial_t u-\Delta u=\varphi(t)|u|^p$, the initial boundary value problem on domains of $\mathbb{R}^n$ (with Dirichlet boundary conditions) is studied in \cite{Mei90}.

 Later, the equation $\partial_t u -\Delta_{\mathrm{\bf X}} u=\varphi(t) f(u)$, where $\Delta_{\mathrm{\bf X}}=\sum_{j}X_j^2$ has been studied for a system of left-invariant vector fields $\mathrm{\bf X}=\{X_j\}_j$ on a unimodular Lie group satisfying the H\"ormander condition in \cite{CKR24} and for a system of vector fields $\mathrm{\bf X}=\{X_j\}_j$ on $\mathbb{R}^n$ which are homogeneous with respect to a suitable family of dilations and satisfying the H\"ormander condition in \cite{CKR26}. 

For the Cauchy problem associated to the semilinear non-autonomous damped wave equation
\begin{align}\label{semilinear CP damped Eucl nonaut}
\begin{cases} \partial_t^2 u-\Delta u+\partial_t u =\varphi(t)| u|^p, &  x\in \mathbb{R}^n, \ t>0,\\
u(0,x)=\varepsilon u_0(x), & x\in \mathbb{R}^n, \\ \partial_t u(0,x)=\varepsilon u_1(x), & x\in \mathbb{R}^n,
\end{cases}
\end{align}
 in the polynomial case $\varphi(t)=(1+t)^\alpha$, with $\alpha>-1$, the critical exponent is given by $p_{\mathrm{Fuj}}(\frac{n}{\alpha+1})$, cf. \cite{D13,DL13}.

In the present work, we focus on the model \eqref{semilinear CP damped phi} on a compact Lie group, where the situation is quite different from the Euclidean one.

\paragraph*{Notations}

In the next sections, we denote by $L^q(\mathbb{G})$ the space of $q$-summable functions on $\mathbb{G}$ with respect to the normalized Haar measure for $1\leqslant q$ (respectively, the space of essentially bounded functions for  $q=\infty$). Moreover, for $s>0$ and $q>1$ the Sobolev space $H^{s,q}_\mathcal{L}(\mathbb{G)}$ is given by $$H^{s,q}_\mathcal{L}(\mathbb{G)}\doteq \left\{ f\in L^q(\mathbb{G}): (-\mathcal{L})^{s/2}f\in L^q(\mathbb{G})\right\}$$ endowed with the norm
$$\|f\|_{H^{s,q}_\mathcal{L}(\mathbb{G)}}\doteq\| f\|_{L^q(\mathbb{G})}+\| (-\mathcal{L})^{s/2} f\|_{L^q(\mathbb{G})}.$$ When $q=2$, we denote the Hilbert space $H^{s,2}_\mathcal{L}(\mathbb{G)}$ simply by $H^{s}_\mathcal{L}(\mathbb{G)}$. \\ Finally, we write $f\lesssim g$ when there exists a positive constant $C$ such that $f\leqslant Cg$ and $f\approx g$ when $g\lesssim f \lesssim g$.

\subsection{The uniform upper scaling condition}

In this section we introduce the key assumption on the factor $\varphi$ that will allow us to prove the blow-up result.

\begin{definition} \label{Def UUSC} We say that the function $\varphi: [0,+\infty)\to [0,+\infty)$ satisfies the \emph{uniform upper scaling condition} if there exist $K>0$,  $\lambda_0 \in (0,1)$,  $a \geqslant 0$ such that $ \varphi(s)\geqslant K \varphi(\lambda s)$ for any $\lambda \in (\lambda_0,1]$ and any $s\geqslant a$.
\end{definition}

Although the uniform upper scaling condition might look rather technical, we may construct several examples of functions $\varphi:[0,+\infty)\to [0,+\infty)$ that satisfy it. For the sake of brevity, in this section we denote $$\mathfrak{U}\doteq \{\varphi:[0,+\infty)\to [0,+\infty)\ |\ \varphi \mbox{ satisfies the uniform upper scaling condition}\}.$$
\begin{enumerate}
\item If $\varphi$ is an increasing function, then $\varphi \in \mathfrak{U}$ (with $K=1$, $a=0$ and for any given $\lambda_0\in(0,1)$).
\item If $\varphi(s)\doteq s^{\alpha}$ for $s\gg 1$ with $\alpha\in \mathbb{R}$, then  $\varphi \in \mathfrak{U}$ (with $K=\lambda_0^{\max\{-\alpha,0\}}$ and given $a>0$ and $\lambda_0 \in (0,1)$).
\item If $\varphi$ is a definitively positive function satisfying 
$$\displaystyle{\limsup_{s\to +\infty} \frac{\varphi(\lambda s)}{\varphi(s)}\in (0,+\infty) \quad \mbox{ uniformly with respect to } \lambda \in (\lambda_0,1]},$$  then $\varphi \in \mathfrak{U}$.
\item If $\varphi(s)\doteq (\ln s)^{\beta}$ for $s\gg 1$ with $\beta \in \mathbb{R}$, then $\varphi \in\mathfrak{U}$. Indeed, if $\beta\geqslant 0$, then $\varphi$ is increasing, while if $\beta<0$, then
\begin{align*}
\frac{\varphi(\lambda s)}{\varphi(s)}=\left(1+\frac{\ln \lambda}{\ln s}\right)^\beta\leqslant \left(1+\frac{\ln \lambda_0}{\ln s}\right)^\beta \longrightarrow 1
\end{align*} as $s\to +\infty$ uniformly with respect to $\lambda \in (\lambda_0,1]$ (here $\lambda_0$ is any given element of $(0,1)$).
\item If $\varphi(s)\doteq  \mathrm{e}^{\gamma s}$ with $\gamma\in \mathbb{R}$, then for $\gamma\geqslant 0$ we have $\varphi \in\mathfrak{U}$ (due to the monotonicity), while for $\gamma<0$ we have $\varphi \not\in\mathfrak{U}$, since 
\begin{align*}
\frac{\varphi(\lambda s)}{\varphi(s)}=\mathrm{e}^{\gamma (\lambda-1)s}\longrightarrow +\infty
\end{align*} as $s\to +\infty$ for any $\lambda\in (0,1)$.
\item Let $\beta \in\mathbb{R}$. We define $\varphi_\beta(s)\doteq (1+s)^{-1} (\ln(\mathrm{e}+s))^{-\beta}$ for any $s\geqslant 0$. Since
\begin{align*}
\frac{\varphi_{\beta}(\lambda s)}{\varphi_{\beta}(s)}=\frac{1+s}{1+\lambda s}\left(\frac{\ln(\mathrm{e}+s)}{\ln(\mathrm{e}+\lambda s)}\right)^{\beta},
\end{align*} proceeding as in item 4 of this list, we have that $\varphi_\beta \in\mathfrak{U}$.
\item Let $\alpha \in\mathbb{R}$. We define $\chi_\alpha(s)\doteq (1+s)^\alpha [2+\sin(s^2)]$  for any $s\geqslant 0$. For any $s\geqslant 0$ we have $2+\sin(s^2)\in [1,3]$. Therefore,
\begin{align*}
\frac{\chi_\alpha(\lambda s)}{\chi_\alpha(s)}\leqslant 3 \left(\frac{1+\lambda s}{1+s}\right)^\alpha,
\end{align*} and hence, reasoning a before, we conclude that $\chi_\alpha\in\mathfrak{U}$. \\ Notice that, differently from all previous examples, $\chi_\alpha$ is not monotone.
\end{enumerate}

\begin{remark}[convex cone property]
If $\varphi_1,\varphi_2\in\mathfrak{U}$ and $\alpha_1,\alpha_2\geqslant 0$, then $\alpha_1\varphi_1+\alpha_2\varphi_2\in \mathfrak{U}$.
\end{remark}

\begin{remark}[algebra property]
If $\varphi_1,\varphi_2\in\mathfrak{U}$, then $\varphi_1\cdot \varphi_2\in \mathfrak{U}$.
\end{remark}

\begin{remark} We point out that in item 3. of the previous list, we are assuming that $\varphi$ is definitely positive. A significant class of functions $\varphi$ that do not fulfill Definition \ref{Def UUSC} contains functions whose zero sets are unbounded. On the other hand, if the zeros are concentrated in a bounded set, the uniform upper scaling condition may be checked on a half-line where the function is positive. 
\end{remark}

\subsection{Main theorems}

Given $\varphi\in L^1_{\mathrm{loc}}([0,+\infty))$, we denote for any $t\geqslant 0$
\begin{align}\label{def Phi}
\Phi(t) \doteq \int_0^t \varphi(s) \,\mathrm{d}s.
\end{align} 
Let us begin by stating the well-posedness results in energy spaces for the semilinear Cauchy problem \eqref{semilinear CP damped phi}.
For the definition of \emph{mild solutions} see Section \ref{subsection Duhamel}.

\begin{theorem}[Global existence]  \label{Thm global esistence}
Let $\mathbb{G}$ be a compact, connected Lie group and let $n$ be the topological dimension of $\mathbb{G}$. Let us assume $n\geqslant 3$. Let $(u_0,u_1)\in H_\mathcal{L}^1(\mathbb{G})\times L^2(\mathbb{G})$. Let $\varphi\in L^1([0,+\infty))$ be a nonnegative function and $p>1$ such that $p\leqslant \frac{n}{n-2}$. \\ Then, there exists $\varepsilon_0=\varepsilon_0(u_0,u_1,p,\varphi)>0$  such that for any $\varepsilon \in(0,\varepsilon_0]$ the Cauchy problem \eqref{semilinear CP damped phi} admits a uniquely determined mild solution $$u\in \mathcal{C}\left([0,+\infty),H^1_\mathcal{L}(\mathbb{G})\right)\cap \mathcal{C}^1\left([0,+\infty),L^2(\mathbb{G})\right).$$
\end{theorem}

\begin{theorem}[Local existence]  \label{Thm loc esistence}
Let $\mathbb{G}$ be a compact, connected Lie group and let $n$ be the topological dimension of $\mathbb{G}$. Let us assume $n\geqslant 3$. Let $(u_0,u_1)\in H_\mathcal{L}^1(\mathbb{G})\times L^2(\mathbb{G})$. Let $\varphi\in L^1_{\mathrm{loc}}([0,+\infty))$ be a nonnegative function such that $\varphi\not \in L^1([0,+\infty))$ and $p>1$ such that $p\leqslant \frac{n}{n-2}$.\\
If $T>0$ satisfies $\Phi(T)\lesssim \varepsilon^{-(p-1)}$, then the Cauchy problem \eqref{semilinear CP damped phi} admits a uniquely determined mild solution $$u\in \mathcal{C}\left([0,T],H^1_\mathcal{L}(\mathbb{G})\right)\cap \mathcal{C}^1\left([0,T],L^2(\mathbb{G})\right).$$ In particular, if $\varphi$ is a positive function, then, the lifespan $T(\varepsilon)$ satisfies the following lower bound estimate:
\begin{align}\label{lower bound estimate}
T(\varepsilon) \geqslant \Phi^{-1}\left(C_1 \varepsilon^{-(p-1)}\right),
\end{align} where the multiplicative constant $C_1>0$ is independent of $\varepsilon$ and $\Phi$ is defined in \eqref{def Phi}.

\end{theorem}


Next, we state the blow-up result. Before doing so, we introduce the notion of energy solutions for the semilinear Cauchy problem \eqref{semilinear CP damped phi} and the definition of uniform upper scaling condition for $\varphi$.

\begin{definition} \label{Definition energy sol}
Let $\varepsilon>0$ and $(u_0,u_1)\in H_\mathcal{L}^1(\mathbb{G})\times L^2(\mathbb{G})$. We say that $$u\in \mathcal{C}\left([0,T),H^1_\mathcal{L}(\mathbb{G})\right)\cap \mathcal{C}^1\left([0,T),L^2(\mathbb{G})\right) \quad \mbox{such that} \quad \varphi |u|^p \in L^1_{\mathrm{loc}}\left([0,T)\times\mathbb{G}\right)$$ is an \emph{energy solution} on $[0,T)$ to \eqref{semilinear CP damped phi} if $u$ fulfills the integral relation
\begin{align}
& \int_{\mathbb{G}} \partial_t u(t,x) \psi(t,x)\, \mathrm{d}x -\int_{\mathbb{G}}  u(t,x) \psi_s(t,x)\, \mathrm{d}x+\int_{\mathbb{G}}  u(t,x) \psi(t,x)\, \mathrm{d}x \notag \\
 & \quad -\varepsilon \int_{\mathbb{G}} u_1(x) \psi(0,x)\, \mathrm{d}x +\varepsilon \int_{\mathbb{G}}  u_0(x) \psi_s(0,x)\, \mathrm{d}x-\varepsilon \int_{\mathbb{G}}  u_0(x) \psi(0,x)\, \mathrm{d}x \notag \\
 & \quad + \int_0^t\int_\mathbb{G} u(s,x) \big(\psi_{ss}(s,x)-\mathcal{L}\psi(s,x)-\psi_s(s,x)\big) \mathrm{d}x \, \mathrm{d}s = \int_0^t\int_\mathbb{G} \varphi(s)|u(s,x)|^p \psi(s,x)\,  \mathrm{d}x \, \mathrm{d}s \label{def energ sol int relation}
\end{align} for any $\psi\in\mathcal{C}^\infty_0([0,T)\times \mathbb{G})$ and any $t\in (0,T)$. \\
Furthermore, we define the \emph{lifespan}  of $u$ as $T(\varepsilon)\doteq \sup \left\{T>0: u \mbox{ is an energy solution to \eqref{semilinear CP damped phi} on } [0,T) \right\}$.
\end{definition}

\begin{theorem}[Blow-up]  \label{Thm blow up}
Let $\mathbb{G}$ be a compact Lie group. Let $(u_0,u_1)\in H_\mathcal{L}^1(\mathbb{G})\times L^2(\mathbb{G})$ satisfying
\begin{align}\label{assumption integral cauchy data}
\int_{\mathbb{G}} u_j(x) \,\mathrm{d}x>0 \qquad \mbox{ for } j=0,1.
\end{align} 
Let $p>1$ and let us assume that $\varphi\in L^1_{\mathrm{loc}}([0,+\infty))$ is a function satisfying the uniform upper scaling condition such that $\varphi \not \in L^1([0,+\infty))$. \\ Let $u$ be an energy solution to \eqref{semilinear CP damped phi} 
with lifespan $T=T(\varepsilon)$.  
Then, $u$ blows up in finite time.
\\ Moreover, if $\varphi$ is a positive function, there exists a positive constant $\varepsilon_0=\varepsilon_0(u_0,u_1,p,\varphi)>0$ such that for any $\varepsilon\in (0,\varepsilon_0]$  the upper bound estimate for the lifespan 
\begin{equation}\label{Lifespan upper bound Thm statement}
T(\varepsilon)\leqslant c_2\, \Phi^{-1} \left(C_2 \varepsilon^{-(p-1)}\right)
\end{equation} holds, where the multiplicative constants $c_2>1$, $C_2>0$ are independent of $\varepsilon$ and $\Phi$ is defined in \eqref{def Phi}.

\end{theorem}

\begin{remark}
The positivity of the function $\varphi$ in Theorem \ref{Thm blow up} is assumed just to write the upper bound estimate for the lifespan in a pleasant way (i.e. by using the inverse function of $\Phi$), but it is completely unnecessary to prove the blow-up.
\end{remark}

\begin{remark}
We point out that the positivity assumption in \eqref{assumption integral cauchy data} can be actually weakened by requiring that $\int_{\mathbb{G}}u_0(x) \mathrm{d}x>0$ and $ \int_{\mathbb{G}}u_1(x) \mathrm{d}x\geqslant 0$.
\end{remark}

\begin{remark} Let $\beta\in \mathbb{R}$ and let for any $s\geqslant 0$
\begin{align*}
\varphi_{\beta}(s)\doteq \frac{1}{(1+s)[\ln(\mathrm{e}+s)]^{\beta}}.
\end{align*}
Clearly $\varphi_{\beta}\in L^1([0,+\infty))$ if and only if $\beta>1$. Therefore, since any mild solution is an energy solution, combining the results from Theorems \ref{Thm global esistence}, \ref{Thm loc esistence} and \ref{Thm blow up}, we obtain the following lifespan estimates:
\begin{align*}
	T_\beta(\varepsilon)\approx\begin{cases}
	\exp\left(C\,\varepsilon^{-\frac{p-1}{1-\beta}}\right)&\mbox{if}\ \ \beta<1,\\[0.5em]
	\exp\left(\exp\left(C\,\varepsilon^{-(p-1)}\right)\right)&\mbox{if}\ \ \beta=1, \\[0.5em]
	+\infty &\mbox{if}\ \ \beta>1.
	\end{cases}
\end{align*}
In particular, the double exponential decay in the limit case $\beta=1$ emphasizes how the summability for $\varphi_\beta$ is a criterion to separate between the blow-up and the global existence of small data solutions to \eqref{semilinear CP damped phi}.
\end{remark}

\section{Existence of mild solutions}
\subsection{Integral formulation and Duhamel's principle}\label{subsection Duhamel}

In the next sections, we prove Theorems \ref{Thm global esistence} and \ref{Thm loc esistence}. Let us get started by recalling the notion of mild solutions to \eqref{semilinear CP damped phi}. By using Duhamel's principle, we represent the solution to the linear inhomogeneous problem
\begin{align}\label{linear CP damped}
\begin{cases} \partial_t^2 u-\mathcal{L} u+\partial_t u =F(t,x), &  x\in \mathbb{G}, \ t>0,\\
u(0,x)= u_0(x), & x\in \mathbb{G}, \\ \partial_t u(0,x)= u_1(x), & x\in \mathbb{G}.
\end{cases}
\end{align} More precisely, denoting by $E_0(t,x)$ and $E_1(t,x)$ the fundamental solutions to \eqref{linear CP damped} in the homogeneous case $F= 0$ with initial data $(u_0,u_1)=(\delta_0,0)$ and $(u_0,u_1)=(0,\delta_0)$, respectively, the solution to \eqref{linear CP damped} is given by
\begin{align*}
u(t,x) = u_0(x)\ast_{(x)} E_0(t,x) + u_1(x)\ast_{(x)} E_1(t,x)+ \int_0^t F(s,x)\ast_{(x)} E_1(t-s,x) \, \mathrm{d}s.
\end{align*}
We emphasize that we employed the invariance by time translations for the differential operator $\partial_t^2-\mathcal{L}+\partial_t$ and the convolution identity $\mathcal{D}\big(v\ast_{(x)}E_1(t,\cdot)\big)=v\ast_{(x)}\mathcal{D}(E_1(t,\cdot))$ for any left-invariant differential operator $\mathcal{D}$ on $\mathbb{G}$.

We say that $u$ is a \emph{mild solution} to \eqref{semilinear CP damped phi} on $[0,T]$ if $u$ is a fixed point for the nonlinear integral operator 
\begin{align*}
 N u(t,x) \doteq \varepsilon u_0(x)\ast_{(x)} E_0(t,x) +\varepsilon
 u_1(x)\ast_{(x)} E_1(t,x)+ \int_0^t \varphi(s) |u(s,x)|^p\ast_{(x)} E_1(t-s,x) \, \mathrm{d}s
\end{align*}
 in the  space $X(T)\doteq\mathcal{C}\left([0,T],H^1_\mathcal{L}(\mathbb{G})\right)\cap \mathcal{C}^1\left([0,T],L^2(\mathbb{G})\right)$, endowed with the following norm:
$$\|u\|_{X(T)}\doteq \sup_{t\in[0,T]} \left(\|u(t,\cdot)\|_{L^2(\mathbb{G})}+\|(-\mathcal{L})^{1/2}u(t,\cdot)\|_{L^2(\mathbb{G})}+\|\partial_t u(t,\cdot)\|_{L^2(\mathbb{G})}\right).$$
In the choice of the norm of $X(T)$ we have a loss of decay rate for the $L^2(\mathbb{G})$ - norm of $(-\mathcal{L})^{1/2}u$ and $\partial_t u$ in comparison to the corresponding homogeneous problem. Indeed we have the following result. 

\begin{proposition} \label{Prop L^2-L^2 estimates}
Let $(u_0,u_1)\in H^1_\mathcal{L}(\mathbb{G})\times L^2(\mathbb{G})$ and let $u\in \mathcal{C}\big([0,\infty),H^1_\mathcal{L}(\mathbb{G})\big)\cap  \mathcal{C}^1\big([0,\infty),L^2(\mathbb{G})\big)$ be the solution to the homogeneous Cauchy problem
\begin{align} \label{linear CP damped hom}
\begin{cases} \partial_t^2 u-\mathcal{L} u+\partial_t u =0, &  x\in \mathbb{G}, \ t>0,\\
u(0,x)= \varepsilon u_0(x), & x\in \mathbb{G}, \\ \partial_t u(0,x)= \varepsilon u_1(x), & x\in \mathbb{G}.
\end{cases}
\end{align} Then, $u$ satisfies the following $L^2(\mathbb{G})$ - $L^2(\mathbb{G})$ estimates
\begin{align*}
\|u(t,\cdot)\|_{L^2(\mathbb{G})} &\leqslant C \varepsilon\left( \|u_0\|_{L^2(\mathbb{G})}+\|u_1\|_{L^2(\mathbb{G})}\right), 
\\
\|(-\mathcal{L})^{1/2}u(t,\cdot)\|_{L^2(\mathbb{G})} &\leqslant C \varepsilon (1+t)^{-\frac{1}{2}}\left( \|u_0\|_{H^1_\mathcal{L}(\mathbb{G})}+\|u_1\|_{L^2(\mathbb{G})}\right), 
\\
\| \partial_t u(t,\cdot)\|_{L^2(\mathbb{G})} &\leqslant C \varepsilon (1+t)^{-1}\left( \|u_0\|_{H^1_\mathcal{L}(\mathbb{G})}+\|u_1\|_{L^2(\mathbb{G})}\right), 
\end{align*} for any $t\geqslant 0$, where $C$ is a positive multiplicative constant independent of $t$ and $\varepsilon$.
\end{proposition}

The proof of the previous proposition can be found in \cite{Pal20DWE}.
In particular, these $L^2(\mathbb{G})$ - $L^2(\mathbb{G})$ estimates for the solution to \eqref{linear CP damped hom} are derived by using the group Fourier transform with respect to the spatial variable $x$. Indeed, thanks to Plancherel identity, it is possible to determine by duality an explicit representation for the $L^2(\mathbb{G})$ norms of $u(t,\cdot)$, $(-\mathcal{L})^{1/2}u(t,\cdot)$ and $\partial_t u(t,\cdot)$. We refer to \cite{RT10} for further details on the representation theory on compact Lie groups. The idea to work with the group Fourier transform generalizes in a natural way the approach in the Euclidean framework (cf. \cite{Mat76}), and it has been introduced for linear wave models on compact Lie groups in \cite{GR15} and later generalized in the case of the damped wave equation to other settings (see, for example, \cite{RT18,Pal19} on the Heisenberg group or \cite{ITW24} on a measure space by means of the spectral resolution for the nonnegative and self-adjoint operator replacing the Laplacian).

In order to show that $N$ admits a uniquely determined fixed point (either for small Cauchy data in Theorem \ref{Thm global esistence} or for a time $T$ below a certain $\varepsilon$-depending condition in Theorem \ref{Thm loc esistence}), we will employ Banach's fixed point theorem. To deal with the nonlinear term in the space $L^2(\mathbb{G})$, we apply a Gagliardo-Nirenberg type inequality derived in \cite{RY19} for the more general framework of connected Lie groups (cf. Lemma \ref{Lemma GN ineq L2}).

\subsection{Proof of Theorem \ref{Thm global esistence}}  \label{Subsection GESDS}

A fundamental tool to prove the global/local existence results is a Gagliardo-Nirenberg type inequality on compact Lie groups. The next lemma is a very special case of \cite[Theorem 1.5]{RY19} (see also Lemma 2.2 and Remark 4 in \cite{Pal20DWE} for further clarifications).

\begin{lemma}\label{Lemma GN ineq L2}
Let $\mathbb{G}$ be a connected unimodular Lie group with topological dimension $n\geqslant 3$. For any $q\geqslant 2$ such that $q\leqslant \frac{2n}{n-2}$ the following Gagliardo-Nirenberg type inequality holds
\begin{align} \label{GN ineq L2}
\| f\|_{L^q(\mathbb{G})}\lesssim \| f\|^{\theta(n,q)}_{H^{1}_\mathcal{L}(\mathbb{G})} \| f\|^{1-\theta(n,q)}_{L^{2}(\mathbb{G})}
\end{align}
for any $f\in H^{1}_\mathcal{L}(\mathbb{G})$, where $ \theta(n,q)\doteq  n\left(\frac{1}{2}-\frac{1}{q}\right)$.
\end{lemma}

Under the assumption $\varphi \in L^1([0,+\infty))$, we now prove the existence of a uniquely determined solution to \eqref{semilinear CP damped phi} in $X(T)$ with estimates that are uniform with respect to $T$, so that the solution can be globally in time prolonged. \\
Let us estimate $\|Nu\|_{X(T)}$ for $u\in X(T)$. We begin by  rewriting $$Nu=  u^{\mathrm{ln}}+ I u,$$ where $$ u^{\mathrm{ln}}(t,x)\doteq \varepsilon u_0(x)\ast_{(x)} E_0(t,x) +\varepsilon  u_1(x)\ast_{(x)} E_1(t,x) $$ is the solution to the linear homogeneous problem \eqref{linear CP damped hom} and $$I u(t,x) \doteq \int_0^t \varphi(s) |u(s,x)|^p\ast_{(x)} E_1(t-s,x) \, \mathrm{d}s$$ is the Duhamel's integral term. \\
 By using Proposition \ref{Prop L^2-L^2 estimates} we find that $$\|u^{\mathrm{ln}}\|_{X(T)}\lesssim  \varepsilon\,  \|(u_0,u_1)\|_{H^1_\mathcal{L}(\mathbb{G})\times L^2(\mathbb{G})}.$$ On the other hand, being the linear Cauchy problem \eqref{linear CP damped hom} invariant by time translations, we obtain
\begin{align}
\|\partial_t^j (-\mathcal{L})^{i/2} Iu(t,\cdot)\|_{L^2(\mathbb{G})}  & \lesssim \int_0^t \varphi(s) (1+t-s)^{-j-\frac{i}{2}} \| u(s,\cdot)\|^p_{L^{2p}(\mathbb{G})} \, \mathrm{d}s \notag \\ 
& \lesssim \int_0^t \varphi(s) \| u(s,\cdot)\|^{p\theta(n,2p)}_{H^{1}_{\mathcal{L}}(\mathbb{G})} \| u(s,\cdot)\|^{p(1-\theta(n,2p))}_{L^2(\mathbb{G})} \, \mathrm{d}s \notag \\ &\lesssim \int_0^t \varphi(s)\, \mathrm{d}s  \, \| u\|_{X(t)}^p \leqslant \|\varphi\|_{L^1([0,+\infty))} \, \| u\|_{X(T)}^p \label{estimate Ju in X(T)}
\end{align} for $i,j\in\{0,1\}$ such that $0\leqslant i+j\leqslant 1$. \\ We stress that the upper bound for $p\leqslant \frac{n}{n-2}$ in the statement of Theorem \ref{Thm loc esistence} is due to the employment of the Gagliardo-Nirenberg inequality from Lemma \ref{Lemma GN ineq L2} in the previous estimate. Analogously, combining the estimate $$||u|^p-|v|^p|\lesssim |u-v|(|u|^{p-1}+|v|^{p-1}),$$ H\"older's inequality and \eqref{GN ineq L2}, we obtain for $i,j\in\{0,1\}$ such that $0\leqslant i+j\leqslant 1$ 
\begin{align}
\|\partial_t^j (-\mathcal{L})^{i/2} (Iu(t,\cdot)-Iv(t,\cdot))\|_{L^2(\mathbb{G})}  & \lesssim \int_0^t \varphi(s) (1+t-s)^{-j-\frac{i}{2}} \| |u(s,\cdot)|^p-|v(s,\cdot)|^p\|_{L^{2}(\mathbb{G})} \, \mathrm{d}s \notag \\  
& \lesssim \int_0^t \varphi(s) \| u(s,\cdot)-v(s,\cdot)\|_{L^{2p}(\mathbb{G})} \left(\|u(s,\cdot)\|^{p-1}_{L^{2p}(\mathbb{G})}+ \|v(s,\cdot)\|^{p-1}_{L^{2p}(\mathbb{G})}\right) \, \mathrm{d}s \notag \\ 
& \lesssim \int_0^t \varphi(s)\, \mathrm{d}s  \, \| u-v\|_{X(t)}\left(\|u\|^{p-1}_{X(t)}+\|v\|^{p-1}_{X(t)}\right) \notag \\ 
& \, \leqslant \|\varphi\|_{L^1([0,+\infty))} \| u-v\|_{X(T)}\left(\|u\|^{p-1}_{X(T)}+\|v\|^{p-1}_{X(T)}\right). \label{estimate Ju -Jv in X(T)}
\end{align} Summarizing, we proved that 
\begin{align*}
\| N u\|_{X(T)} &\leqslant C  \varepsilon\,  \|(u_0,u_1)\|_{H^1_\mathcal{L}(\mathbb{G})\times L^2(\mathbb{G})} +C \|\varphi\|_{L^1([0,+\infty))} \|u\|^p_{X(T)}, \\
\| N u -Nv\|_{X(T)} &\leqslant C \|\varphi\|_{L^1([0,+\infty))} \| u-v\|_{X(T)}\left(\|u\|^{p-1}_{X(T)}+\|v\|^{p-1}_{X(T)}\right).
\end{align*} Therefore, for $\varepsilon$ sufficiently small and depending on $\|\varphi\|_{L^1([0,+\infty))}$, $N$ is a contraction on a certain ball around $0$ in the Banach space $X(T)$, so Banach's fixed point provides a uniquely determined fixed point $u$ for $N$ which is exactly our mild solution to \eqref{semilinear CP damped phi} on $[0,T]$. Since the multiplicative constant $C$ is independent of $T$, our mild solution $u$ can be prolonged to $[0,+\infty)$.

\subsection{Proof of Theorem \ref{Thm loc esistence}} 

In this section, we prove the local existence result for $\varphi\in L^1_{\mathrm{loc}}([0,+\infty))$ such that $\varphi \not \in L^1([0,+\infty))$.

We can repeat the same considerations from Section \ref{Subsection GESDS}, with a very minor modification in \eqref{estimate Ju in X(T)} and in \eqref{estimate Ju -Jv in X(T)}: since we are now working under the assumption that $\varphi$ is not summable, by using the function $\Phi$ defined in \eqref{def Phi}, for $i,j\in\{0,1\}$ such that $0\leqslant i+j\leqslant 1$ we obtain
\begin{align*}
 \|\partial_t^j (-\mathcal{L})^{i/2} Iu(t,\cdot)\|_{L^2(\mathbb{G})}   & \lesssim \Phi(t)  \, \| u\|_{X(t)}^p, \\
\|\partial_t^j (-\mathcal{L})^{i/2} (Iu(t,\cdot)-Iv(t,\cdot))\|_{L^2(\mathbb{G})}  & \lesssim \Phi(t)  \,  \| u-v\|_{X(t)}\left(\|u\|^{p-1}_{X(t)}+\|v\|^{p-1}_{X(t)}\right).
\end{align*} Hence, for a suitable positive multiplicative constant $C$ we have
\begin{align*}
\| N u\|_{X(t)} &\leqslant C  \varepsilon\,  \|(u_0,u_1)\|_{H^1_\mathcal{L}(\mathbb{G})\times L^2(\mathbb{G})} +C \Phi(t) \, \|u\|^p_{X(t)}, \\
\| N u -Nv\|_{X(t)} &\leqslant C \Phi(t) \| u-v\|_{X(t)}\left(\|u\|^{p-1}_{X(t)}+\|v\|^{p-1}_{X(t)}\right).
\end{align*}
Setting $R_0\doteq \|(u_0,u_1)\|_{H^1_\mathcal{L}(\mathbb{G})\times L^2(\mathbb{G})}$, we define $R \doteq 2CR_0$ and $\mathfrak{B}(R\varepsilon,t)\doteq \{u\in X(t): \|u\|_{X(t)}\leqslant R\varepsilon \}$. \\  Then, for any $t>0$ such that $\Phi(t)\leqslant  \varepsilon^{-(p-1)} /(4CR^{p-1})$ and any $u,v\in \mathfrak{B}(R\varepsilon,t)$
 and, we have $\|Nu\|_{X(t)}\leqslant R\varepsilon$ and $\| N u -Nv\|_{X(t)}\leqslant  \frac{1}{2} \|u-v\|_{X(t)}$. In other words, for any $t>0$ such that $\Phi(t)\leqslant \varepsilon^{-(p-1)} /(4CR^{p-1})$ the operator $N$ is a contraction mapping on $\mathfrak{B}(R\varepsilon,t)$. This completes the local existence and the lower bound estimate for the lifespan.

\section{Blow-up result}

Next, we prove Theorem \ref{Thm blow up} by using an iteration argument. In particular, we will adapt to our model the slicing procedure for iterative schemes introduced in \cite{AKT00} and adapted in \cite{ChenPal19MGT} to handle unbounded exponential multipliers. The peculiarity of our approach relies on the fact that while performing the slicing procedure we lower the upper bound of the domain of integration as well.

Let $u$ be a local in time energy solution to \eqref{semilinear CP damped phi} according to Definition \ref{Definition energy sol} with lifespan $T$. We fix $t\in (0,T)$ and we consider a test function $\psi\in \mathcal{C}^\infty_0 ([0,T)\times \mathbb{G})$ such that $\psi=1$ on $[0,t]\times\mathbb{G}$ in \eqref{def energ sol int relation}. Then,
\begin{align*}
& \int_{\mathbb{G}} \partial_t u(t,x) \, \mathrm{d}x +\int_{\mathbb{G}}  u(t,x) \, \mathrm{d}x  -\varepsilon \int_{\mathbb{G}} u_1(x) \, \mathrm{d}x -\varepsilon \int_{\mathbb{G}}  u_0(x) \, \mathrm{d}x  = \int_0^t \varphi(s)\int_\mathbb{G} |u(s,x)|^p   \mathrm{d}x \, \mathrm{d}s. 
\end{align*} Let us consider the space-average of $u$
\begin{align*}
U_0(t)\doteq \int_{\mathbb{G}}  u(t,x) \, \mathrm{d}x.
\end{align*} 
The evolution of the time-dependent functional $U_0$ will be used to derive the blow-up result.
The previous integral identity can be rewritten as
\begin{align*}
U_0'(t)+U_0(t) - U_0'(0)-U_0(0)=\int_0^t\int_\mathbb{G}  \varphi(s)|u(s,x)|^p   \mathrm{d}x \, \mathrm{d}s \geqslant \int_0^t  \varphi(s) |U_0(s)|^p\, \mathrm{d}s,
\end{align*}  where we applied Jensen's inequality and the fact that the Haar measure on $\mathbb{G}$ is normalized in the last inequality. 
 Multiplying the last relation by $\mathrm{e}^{t}$, we find
\begin{align*}
\frac{\mathrm{d}}{\mathrm{d}t} (\mathrm{e}^{t} U_0(t)) = \mathrm{e}^{t} (U_0'(t)+U_0(t) )\geqslant (U_0'(0)+U_0(0))\, \mathrm{e}^{t}  +\mathrm{e}^{t} \int_0^t  \varphi(s) |U_0(s)|^p\, \mathrm{d}s.
\end{align*} Performing an integration over $[0,t]$, we have
\begin{align*}
\mathrm{e}^{t} U_0(t)\geqslant U_0(0) + (U_0'(0)+U_0(0))\, (\mathrm{e}^{t}  -1)+\int_0^t \mathrm{e}^{\tau} \int_0^\tau  \varphi(s)|U_0(s)|^p\, \mathrm{d}s \, \mathrm{d}\tau.
\end{align*} Summarizing,
\begin{align*} 
 U_0(t)\geqslant U_0(0)+ U_0'(0) (1-\mathrm{e}^{-t} )+\int_0^t \mathrm{e}^{\tau-t} \int_0^\tau  \varphi(s)|U_0(s)|^p\, \mathrm{d}s \, \mathrm{d}\tau.
\end{align*}
Applying Fubini-Tonelli theorem in the double integral, we get
\begin{align*}
\int_0^t \mathrm{e}^{\tau-t} \int_0^\tau  \varphi(s)|U_0(s)|^p\, \mathrm{d}s \, \mathrm{d}\tau & = \mathrm{e}^{-t} \int_0^t \int_s^t \mathrm{e}^{\tau} \varphi(s)|U_0(s)|^p \, \mathrm{d}\tau \, \mathrm{d}s \\ & = \mathrm{e}^{-t} \int_0^t  \varphi(s)|U_0(s)|^p ( \mathrm{e}^{t}-\mathrm{e}^{s}) \, \mathrm{d}s,
\end{align*} and, consequently, for any $t\in[0,T)$
 \begin{align} \label{fundamental ineq U0}
 U_0(t)\geqslant U_0(0)+ U_0'(0) (1-\mathrm{e}^{-t} )+\mathrm{e}^{-t} \int_0^t  \varphi(s) ( \mathrm{e}^{t}-\mathrm{e}^{s})|U_0(s)|^p  \, \mathrm{d}s.
\end{align}
 Indeed, since the non-autonomous nonlinear term is nonnegative, from \eqref{fundamental ineq U0} we have
\begin{align} \label{first lower bound U0}
U_0(t)\geqslant  U_0(0)+ U_0'(0) (1-\mathrm{e}^{-t} ) \geqslant C \varepsilon
\end{align}  for $t\in [0,T)$, where  $C>0$ depends on $u_0,u_1$.  We underline that the assumptions in \eqref{assumption integral cauchy data} guarantee that $C$ can be chosen positive. In particular, $U_0(t)\geqslant 0$ for any $t\in[0,T)$.  \\
Besides, \eqref{fundamental ineq U0} provides the iteration frame for $U_0$, namely, for any $t\in[0,T)$
\begin{align} \label{iteration frame}
U_0(t)\geqslant  \mathrm{e}^{-t} \int_0^t  \varphi(s) ( \mathrm{e}^{t}-\mathrm{e}^{s})(U_0(s))^p  \, \mathrm{d}s.
\end{align}

\subsection{Iteration argument}
In the previous section, we established the iteration frame \eqref{iteration frame} and the first lower bound estimate for $U_0$ in \eqref{first lower bound U0}.
 The next step is to determine a sequence of lower bounds estimates for $U_0$ by using \eqref{iteration frame} in an iterative way. \\
Let us  start by introducing the sequence of positive real numbers $\{\ell_j\}_{j\in\mathbb{N}}$ that characterizes the slicing procedure.
We set $$\ell_j\doteq 1+q^{-j}$$ for any $j\in\mathbb{N}^*$, where $q>1$ is a real parameter. We remark that the infinite product $\prod_{j=1}^{+\infty}\ell_j$ is convergent, being the series $\sum_{j=1}^{+\infty} \ln(1+q^{-j})$ convergent. \\
By using the inequality $\ln(1+\sigma)\leqslant \sigma$ for any $\sigma>-1$ (which follows immediately by the concavity of the logarithmic function) and the sum for the geometric series, we obtain
\begin{align*}
\sum_{j=1}^{+\infty}\ln\left(1+q^{-j}\right)\leqslant \sum_{j=1}^{+\infty}q^{-j} 
=\frac{1}{q-1}.
\end{align*} Hence,
\begin{align*}
\prod_{j=1}^{+\infty}\ell_j = \exp\bigg(\,\sum_{j=1}^{+\infty} \ln\left(1+q^{-j}\right)\bigg)\leqslant \mathrm{e}^{\frac{1}{q-1}}.
\end{align*} 
Let $K>0$, $\lambda_0\in (0,1)$ and $a\geqslant 0$ be the quantities related to $\varphi$ from Definition \ref{Def UUSC}.\\
Since $\displaystyle{\lim_{q\to +\infty} \mathrm{e}^{\frac{1}{q-1}}=1}$ and $\frac{1}{\sqrt{\lambda_0}}>1$, we can choose $q$  sufficiently large so that $\prod_{j=1}^{+\infty}\ell_j\leqslant \frac{1}{\sqrt{\lambda_0}}$. \\ Finally, we pick $$\ell_0\doteq \frac{1}{\sqrt{\lambda_0 }}.$$ \\
In order to simplify the notations, we introduce the sequence of the partial products $\{L_j\}_{j\in\mathbb{N}}$ such that 
\begin{align}\label{Lj def}
L_j\doteq \prod_{k=0}^{j}\ell_k \qquad \mbox{ for any }  j\in\mathbb{N}.
\end{align}
Furthermore, we denote 
\begin{align}\label{L def}
L\doteq \prod_{k=0}^{+\infty}\ell_k \equiv \lim_{j\to +\infty} L_j.
\end{align}
 We point out that $L_j\uparrow L$ as $j\to +\infty$, being $\ell_j>1$ for any $j\geqslant 1$. Thanks to the previous considerations on the choice of $q$, we have $L\leqslant \frac{1}{\lambda_0}$. Due to the monotonicity of the sequence $\{L_j\}_{j\in\mathbb{N}}$, we also have 
\begin{align} \label{Lj upper bound}
L_j \leqslant \frac{1}{\lambda_0} \qquad \mbox{ for any }  j\in\mathbb{N}.
\end{align}

 We stress that the choice of the sequence $\{\ell_j\}_{j\in\mathbb{N}}$ is related to the parameter $\lambda_0$ appearing in the uniform upper scaling condition for $\varphi$ because we are going to use this condition to deal with the $\varphi$ function in the iteration argument. \\
Let us set $$\Phi_a(t) \doteq \int_a^t \varphi(s) \, \mathrm{d}s.$$ 
  Without loss of generality, we may assume $a>0$. \\  We can now write the sequence of lower bound estimates for $U_0$. For any $j\in\mathbb{N}^*$ we have
\begin{align}\label{sequence of lower bound U0}
U_0(t) \geqslant C_j\, \varepsilon^{p^j} \left(\Phi_a\left(\frac{t}{L_{j-1}}\right)\right)^{\frac{p^j-1}{p-1}}  \qquad \mbox{for any} \ t\in [a L_{j-1},T),
\end{align} where $\{C_j\}_{j\geqslant 1}$ is a sequence of nonnegative real numbers that will be determined iteratively in this section. 

We are going to prove \eqref{sequence of lower bound U0} by induction on $j\in\mathbb{N}^*$. Let us begin with the base case $j=1$. 

By plugging \eqref{first lower bound U0} into \eqref{iteration frame}, for $t\in[a \ell_0,T)$ we get
\begin{align*}
U_0(t) & \geqslant C^p \varepsilon^p  \mathrm{e}^{-t} \int_0^t  \varphi(s) ( \mathrm{e}^{t}-\mathrm{e}^{s})\, \mathrm{d}s \\
& \geqslant C^p \varepsilon^p  \mathrm{e}^{-t} \int_a^{\tfrac{t}{\ell_0}}  \varphi(s) ( \mathrm{e}^{t}-\mathrm{e}^{s})  \, \mathrm{d}s \\
& \geqslant C^p \varepsilon^p  \mathrm{e}^{-t} ( \mathrm{e}^{t}-\mathrm{e}^{\tfrac{t}{\ell_0}}) \int_a^{\tfrac{t}{\ell_0}}  \varphi(s)   \, \mathrm{d}s \\
& = C^p \varepsilon^p   \left( 1-\mathrm{e}^{-\tfrac{\ell_0-1}{\ell_0} t}\right) \Phi_a\left(\frac{t}{\ell_0}\right) \geqslant   C^p \left( 1-\mathrm{e}^{-(\ell_0-1)a}\right) \varepsilon^p  \,  \Phi_a\left(\frac{t}{\ell_0}\right). 
\end{align*} Hence, we proved \eqref{sequence of lower bound U0} for $j=1$, with $C_1\doteq  C^p \left( 1-\mathrm{e}^{-(\ell_0-1)a}\right)$.

Next, we prove the inductive step. We assume \eqref{sequence of lower bound U0} satisfied for some $j\geqslant 1$ and we prove it for $j+1$.
Using \eqref{sequence of lower bound U0} in \eqref{iteration frame}, for $t\in [aL_j,T)$ we have
\begin{align*}
U_0(t) & \geqslant  \mathrm{e}^{-t} \int_{aL_{j-1}}^t  \varphi(s) ( \mathrm{e}^{t}-\mathrm{e}^{s})\, (U_0(s))^p \mathrm{d}s \\
& \geqslant C_j^p \varepsilon^{p^{j+1}} \mathrm{e}^{-t} \int_{aL_{j-1}}^t  \varphi(s) ( \mathrm{e}^{t}-\mathrm{e}^{s})\, \left(\Phi_a\left(\frac{s}{L_{j-1}}\right)\right)^{\frac{p^{j}-1}{p-1}p} \mathrm{d}s  \\
& \geqslant C_j^p \varepsilon^{p^{j+1}} \mathrm{e}^{-t} \int_{aL_{j-1}}^{\tfrac{t}{\ell_j}}  \varphi(s) ( \mathrm{e}^{t}-\mathrm{e}^{s})\, \left(\Phi_a\left(\frac{s}{L_{j-1}}\right)\right)^{\frac{p^{j}-1}{p-1}p} \mathrm{d}s  \\
& \geqslant C_j^p \varepsilon^{p^{j+1}} \bigg(1-\mathrm{e}^{-\tfrac{\ell_j-1}{\ell_j}t}\bigg) \int_{aL_{j-1}}^{\tfrac{t}{\ell_j}}  \varphi(s) \, \left(\Phi_a\left(\frac{s}{L_{j-1}}\right)\right)^{\frac{p^{j}-1}{p-1}p} \mathrm{d}s,
\end{align*} where we used  $\{\ell_j\}_{j\in\mathbb{N}}\subset (1,+\infty)$ to lower the endpoint of the interval of integration and the inequality $t\geqslant a L_j$ to guarantee that after the shrinking of the interval of integration from $[aL_{j-1},t]$ to $[aL_{j-1},t/\ell_j]$ the latter interval is not empty.
By using \eqref{Lj upper bound}, we can apply the uniform upper scaling condition for $\varphi$ on $[aL_{j-1},t/\ell_j]$ with $\lambda=1/L_{j-1}$, obtaining for $t\in [aL_j,T)$ 
\begin{align*}
U_0(t) & \geqslant K C_j^p \varepsilon^{p^{j+1}} \bigg(1-\mathrm{e}^{-\tfrac{\ell_j-1}{\ell_j}aL_j}\bigg) \int_{aL_{j-1}}^{\tfrac{t}{\ell_j}}  \varphi\left(\frac{s}{L_{j-1}}\right) \, \left(\Phi_a\left(\frac{s}{L_{j-1}}\right)\right)^{\frac{p^{j}-1}{p-1}p} \mathrm{d}s \\
& \geqslant K C_j^p L_{j-1}\, \varepsilon^{p^{j+1}} \left(1-\mathrm{e}^{-a(\ell_j-1)}\right) \int_{a}^{\tfrac{t}{L_j}}  \varphi(s) \, \left(\Phi_a(s)\right)^{\frac{p^{j}-1}{p-1}p} \mathrm{d}s  \\
& \geqslant K  C_j^p \, \frac{p-1}{p^{j+1}-1}  \left(1-\mathrm{e}^{-a(\ell_j-1)}\right) \varepsilon^{p^{j+1}} \left(\Phi_a\left(\frac{t}{L_j}\right)\right)^{\frac{p^{j+1}-1}{p-1}}, 
\end{align*} where we used the monotonicity of the exponential factor in the first two estimates and $\Phi_a'=\varphi$, $\Phi_a(a)=0$ in the  last inequality.

We notice that 
\begin{align*}
1-\mathrm{e}^{-a(\ell_j-1)}& \geqslant a(\ell_j-1)\left(1-\tfrac{a}{2}(\ell_j-1)\right)\\ & = a q^{-2j}\left(q^j-\tfrac{a}{2}\right) \\ & \geqslant  a \left(q-\tfrac{a}{2}\right) q^{-2j},
\end{align*} where we used the estimate $1-\mathrm{e}^{-\sigma}\geqslant \sigma\left(1-\frac{\sigma}{2}\right)$ for $\sigma\geqslant 0$,  that follows from the Taylor expansion of the exponential function in a right neighborhood of 0. Furthermore, we can always assume (without loss of generality) that $2q>a$ to guarantee that this lower bound for $1-\mathrm{e}^{-a(\ell_j-1)}$ is positive.

Summarizing, for $t\in [aL_j,T)$ we have
\begin{align*}
U_0(t) & \geqslant  \tfrac{K}{2} (p-1) a(2q-a)q^2   (pq^2)^{-(j+1)} C_j^p   \varepsilon^{p^{j+1}} \left(\Phi_a\left(\frac{t}{L_j}\right)\right)^{\frac{p^{j+1}-1}{p-1}}, 
\end{align*} which is exactly \eqref{sequence of lower bound U0} for $j+1$ provided that $C_{j+1}\doteq \frac{K}{2}(p-1) a(2q-a)q^2 (pq^2)^{-(j+1)} C_j^p$.

\subsection{Upper bound estimate for the lifespan}

In this section, we derive the blow-up of $U_0$ from \eqref{sequence of lower bound U0}. We proceed as follows: we show that the $j$-dependent lower bound in \eqref{sequence of lower bound U0} for $U_0$ blows up  as $j\to \infty$ for $t$ greater than a certain $\varepsilon$-dependent threshold. This will also provide, as a byproduct, the upper bound estimate for the lifespan. We begin by deriving a lower bound estimate for multiplicative constant $C_j$ which can be handled more easily.

If we denote $D\doteq  \frac{K}{2}(p-1) a(2q-a)q^2$, we get $C_j= D(pq^2)^{-j}C_{j-1}^p$. Applying the logarithmic function to both sides of the previous relation and using iteratively this equality, we obtain
\begin{align*}
\ln C_j & = p \ln C_{j-1} - j \ln (pq^2) +\ln D \\
& = p^2 \ln C_{j-2} -( j+(j-1)p ) \ln (pq^2) + (1+p)\ln D \\
& =\cdots = p^{j-1}\ln C_1 - \left(\sum_{k=0}^{j-2} (j-k)p^k\right) \ln (pq^2)+ \left(\sum_{k=0}^{j-2} p^k\right)\ln D.
\end{align*} Employing the identities
\begin{align} \label{identity sum (j-k)p^k}
\sum_{k=0}^{j-2} (j-k)p^k = \frac{1}{p-1} \left(\frac{2p-1}{p-1}p^{j-1}-\frac{1}{p-1}-(j+1)\right)\quad\text{and}\quad  \sum_{k=0}^{j-2} p^k = \frac{p^{j-1}-1}{p-1},
\end{align} it follows that
\begin{align*}
\ln C_j = p^{j-1} \left(\ln C_1-\frac{(2p-1) \ln (pq^2)}{(p-1)^2} +\frac{\ln D}{p-1} \right)+\frac{(j+1) \ln (pq^2)}{p-1}+\frac{ \ln (pq^2)}{(p-1)^2}-\frac{\ln D}{p-1}
\end{align*} for any $j\in\mathbb{N} \!\smallsetminus\!\{0,1\}$. \\ 
Let $j_0=j_0(p,q,a,K,L)\in\mathbb{N} \!\smallsetminus\!\{0,1\}$ be the smallest positive integer such that
\begin{align*}
j_0\geqslant \frac{\ln D}{\ln (pq^2)}-\frac{p}{p-1}.
\end{align*}
Therefore, for $j\geqslant j_0$ we find
\begin{align}
\ln C_j &\geqslant p^{j-1} \left(\ln C_1-\frac{(2p-1) \ln (pq^2)}{(p-1)^2} +\frac{\ln D}{p-1} \right) \notag \\& = p^{j-1} \ln \left(D^{1/(p-1)}(pq^2)^{(1-2p)/(p-1)^2}C_1\right)  =p^{j-1} \ln E, \label{lower bound log Cj}
\end{align} where $E\doteq D^{1/(p-1)}(pq^2)^{(1-2p)/(p-1)^2}C_1>0$. \\
Since $L_j \uparrow L$ as $j\to +\infty$,  \eqref{sequence of lower bound U0} holds for any $j\in\mathbb{N}^*$ and any $t\in [ aL,T)$.
Therefore, combining \eqref{sequence of lower bound U0} and \eqref{lower bound log Cj}, we arrive at
\begin{align}
U_0(t) &  \geqslant \exp\left\{p^{j-1}\ln E\right\}\, \varepsilon^{p^j} \left(\Phi_a\left(\frac{t}{L}\right)\right)^{\frac{p^j-1}{p-1}} \notag\\
& = \exp\left\{p^{j}\left[\frac{1}{p}\ln E+\ln \varepsilon +\frac{1}{p-1} \ln \left(\Phi_a\left(\frac{t}{L}\right)\right)\right] \right\} \left(\Phi_a\left(\frac{t}{L}\right)\right)^{-\frac{1}{p-1}} \notag \\
& = \exp\left\{p^{j}\ln \left(H\varepsilon \left(\Phi_a\left(\frac{t}{L}\right)\right)^{\frac{1}{p-1}}\right) \right\} \left(\Phi_a\left(\frac{t}{L}\right)\right)^{-\frac{1}{p-1}}  \label{final lower bound U0}
\end{align}
for any $j\geqslant j_0$ and any $t\in [aL,T)$, where $H\doteq E^{1/p}$. 
Let us introduce the function 
\begin{align*}
J_a(t,\varepsilon)\doteq \ln \left(H\varepsilon \left(\Phi_a\left(t/L\right)\right)^{\frac{1}{p-1}}\right). 
\end{align*} If $J_a(t,\varepsilon)>0$, then, taking the limit as $j\to +\infty$ in \eqref{final lower bound U0} we conclude that $U_0(t)$ cannot be finite. We remark that, given $\varepsilon>0$, we can always find a $t$ sufficiently large such that $J_a(t,\varepsilon)>0$, thanks to the fact that, by assumption, $\varphi$ is not a $L^1([0,+\infty))$ - function and, hence, $\lim_{t\to +\infty}\Phi_a(\tfrac{t}{L})=+\infty$.

Finally, we determine the upper bound estimate for the lifespan.
We recall the notation $\Phi=\Phi_0$ for the primitive of $\varphi$ defined in \eqref{def Phi}. \\
Since $\Phi:[0,+\infty)\to [0,+\infty)$ is a strictly increasing bijection, being $\Phi'=\varphi>0$ and $\varphi \not \in L^1([0,+\infty))$, it follows that $\Phi^{-1}$ is strictly increasing and that $\lim_{\sigma\to +\infty}\Phi^{-1}(\sigma)=+\infty$. \\
Therefore, we can fix a sufficiently small $\varepsilon_0=\varepsilon_0(u_0,u_1,p,a,K,\lambda_0,q)>0$ such that $$\Phi^{-1}\big((H\varepsilon_0)^{-(p-1)}\big)\geqslant a \quad \Leftrightarrow \quad (H\varepsilon_0)^{-(p-1)}\geqslant \int_0^a \varphi(s) \, \mathrm{d}s.$$ Then, for any $\varepsilon\in(0,\varepsilon_0]$ and any $t>L \Phi^{-1}\left(2(H\varepsilon)^{-(p-1)}\right)$, we have
\begin{align*}
t> L \Phi^{-1}\big(2(H\varepsilon)^{-(p-1)}\big) \geqslant L \Phi^{-1}\big(2(H\varepsilon_0)^{-(p-1)}\big) \geqslant La
\end{align*} and, since $\Phi(\frac{t}{L})> 2(H\varepsilon)^{-(p-1)}$,
\begin{align*}
\Phi_a\left(\tfrac{t}{L}\right)=\Phi\left(\tfrac{t}{L}\right) -\int_0^a \varphi(s)\, \mathrm{d}s  & \geqslant \Phi\left(\tfrac{t}{L}\right) -(H\varepsilon_0)^{-(p-1)} \\ & \geqslant \Phi\left(\tfrac{t}{L}\right) -(H\varepsilon)^{-(p-1)} > (H\varepsilon)^{-(p-1)} 
\end{align*} which implies $J_a(t,\varepsilon)>0$ and, therefore, letting $j\to +\infty$ in \eqref{final lower bound U0} we conclude the validity of \eqref{Lifespan upper bound Thm statement} with $c_2=L$ and $C_2=2H^{1-p}$.

\section{Final remarks and open problem}

Combining the results from Theorems \ref{Thm loc esistence} and \ref{Thm blow up}, we see that for a positive not summable $\varphi$, we have the sharp lifespan estimates
\begin{align*}
c_1 \Phi^{-1}\left(C_1 \varepsilon^{-(p-1)}\right) \leqslant T(\varepsilon) \leqslant c_2 \Phi^{-1}\left(C_2 \varepsilon^{-(p-1)}\right)
\end{align*} for any $\varepsilon \in(0,\varepsilon_0]$, where $c_1,c_2,C_1,C_2$ are positive and independent of $\varepsilon$ {constants}. In particular, in the case $\varphi\equiv 1$, our result generalizes the one in \cite{Pal20DWE}.

In the present work, we focus on the model \eqref{semilinear CP damped phi} on a compact Lie group. The situation in this case is relatively simple due to the fact that the critical exponent for the autonomous problem, i.e. for $\varphi\equiv 1$, is $+\infty$.

As we mentioned in the introduction, for the non-autonomous problem \eqref{semilinear CP damped Eucl nonaut} when $\varphi(t)=(1+t)^{\alpha}$ the blow-up has been established in \cite[Example 1.1]{DL13}, while the global existence of small data energy solutions is proved in \cite{D13} (without lifespan estimates).
 To the best of the authors' knowledge, the case of  growing/decaying logarithmic time-dependent factor $\varphi$ remains an open problem. 
A particularly significant case for \eqref{semilinear CP damped Eucl nonaut} is when $p=p_{\mathrm{Fuj}}(n)$ and $\varphi(t)=(\ln (\mathrm{e}+t))^{\beta}$, with $\beta \in\mathbb{R}$. Furthermore, no analogous result seems to be available in the framework of a general stratified Lie group.

Finally, we observe that our results for the damped wave model in \eqref{semilinear CP damped phi} can be easily obtained also for the semilinear heat equation $u_t-\mathcal{L}u=\varphi(t)|u|^p$. For the well-posedness in the energy space, the $L^2(\mathbb{G}) - L^2(\mathbb{G})$ estimates can be obtained either via the group Fourier transform (analogously to what is done in \cite{Pal20DWE}) or directly through the estimates for the heat kernel (cf. \cite[Chapter IV]{Var92}). For the blow-up result, the iteration argument is simpler than in our model. Indeed, for this heat equation the iteration frame is given by $$U_0(t)\geqslant U_0(0)+ \int_0^t \varphi(s) (U_0(s))^p \, \mathrm{d}s.$$ From the above inequality, we understand that no slicing procedure is necessary (as we do not need to deal with an exponential factor in the iteration frame) and, consequently, no uniform upper scaling condition has to be required.

\section*{Acknowledgments}

W. Chen is supported in part by the National Natural Science Foundation of China (grant No. 12301270), Guangdong Basic and Applied Basic Research Foundation (grant No. 2025A1515010240). \\
S. Lucente and A. Palmieri are partially supported by the PRIN 2022 project “Anomalies in partial differential equations and applications” CUP H53C24000820006. \\
S. Lucente and A. Palmieri are members of the Gruppo Nazionale per L’Analisi Matematica, la Probabilità e le loro Applicazioni (GNAMPA) of the Instituto Nazionale di Alta Matematica (INdAM).

\end{document}